\documentclass[twoside,leqno,12pt]{article}
\usepackage{latexsym}
\usepackage{bbm}
\usepackage{amssymb,amsmath}
\begin{document}
\newtheorem{theorem}{Theorem}
\newtheorem{lemma}{Lemma}
\newtheorem{prop}{Proposition}
\newtheorem{corollary}{Corollary}
\newtheorem{conjecture}{Conjecture}
\numberwithin{equation}{section}
\newcommand{\dif}{\mathrm{d}}
\newcommand{\intz}{\mathbb{Z}}
\newcommand{\ratq}{\mathbb{Q}}
\newcommand{\natn}{\mathbb{N}}
\newcommand{\comc}{\mathbb{C}}
\newcommand{\rear}{\mathbb{R}} 
\newcommand{\prip}{\mathbb{P}}
\newcommand{\uph}{\mathbb{H}}
\newcommand{\fie}{\mathbb{F}}

\title{The Lang-Trotter Conjecture on Average}
\date{\today}
\author{Stephan Baier}
\maketitle

\begin{abstract} For an elliptic curve $E$ over $\ratq$ 
and an integer $r$ let 
$\pi_E^r(x)$ be the number of primes $p\le x$ of 
good reduction such that
the trace of the Frobenius morphism of $E/\fie_p$
equals $r$. We consider the quantity $\pi_E^r(x)$
on average over certain sets of elliptic curves. 
More in particular,
we establish the following: If $A,B>x^{1/2+\varepsilon}$ and
$AB>x^{3/2+\varepsilon}$, then the arithmetic
mean of $\pi_E^r(x)$
over all elliptic curves 
$E$ : $y^2=x^3+ax+b$
with $a,b\in \intz$, $|a|\le A$ and $|b|\le B$ is 
$\sim 
C_r\sqrt{x}/\log x$, where $C_r$ is
some constant depending on $r$. This improves a result
of C. David and F. Pappalardi. Moreover, we establish an 
``almost-all'' result on $\pi_E^r(x)$. 
\end{abstract}



\section{Introduction and main results}
Let $E$ be an elliptic curve over $\ratq$. For any prime number $p$ of 
good reduction, let $a_p(E)$ be the trace of the Frobenius morphism of 
$E/\fie_p$. Then the number of points on the reduced curve modulo $p$ 
equals $\#E(\fie_p)= p+1- a_p(E)$. Furthermore, by Hasse's theorem,
$|a_p(E)|\le 2\sqrt{p}$. 

For a fixed integer $r$, let 
$$
\pi_E^r(x):= \#\{p\leq x: a_p(E)=r\}.
$$
If $r=0$ and $E$ has complex multiplication, Deuring \cite{Deu} showed that 

$$
\pi_E^0(x)\sim \frac{\pi(x)}{2} \ \ \ \ \mbox{ as } x\rightarrow\infty.
$$
Primes $p$ with $a_p=0$ are known as ``supersingular primes''.
 
Lang and Trotter \cite{LTr} conjectured that for all other cases an asymptotic 
estimate of the form
 
$$
\pi_E^r(x)\sim C_{E,r}\cdot \frac{\sqrt{x}}{\log x} \ \ \ \ 
\mbox{ as } x\to\infty
$$ 
with a well-defined constant $C_{E,r}\ge 0$ holds. They used a probabilistic 
model to 
give an explicit description of the constant $C_{E,r}$. The constant can be
$0$, and the asymptotic estimate is then interpreted to mean that there is
only a finite number of primes such that $a_p(E)=r$. A concise account of 
Lang-Trotter's probabilistic model and an expression of $C_{E,r}$ as an 
Euler product can be found in \cite{DaP}. 

Fouvry and Murty \cite{FMu} 
obtained average estimates related to the Lang-Trotter
conjecture for the supersingular case $r=0$.
Their result was later generalized by David and Pappalardi \cite{DaP}
to any $r\in \intz$.
In this paper, we shall
improve the results of David and Pappalardi.

As in \cite{DaP}, we define
$$
\pi_{1/2}(x):=\int\limits_2^{x} \frac{{\rm d} t}{2\sqrt{t}\log t}\sim
\frac{\sqrt{x}}{\log x}
$$
and a constant $C_r$ by
\begin{equation}\label{constant}
C_r:=\frac{2}{\pi}\prod\limits_{l|r} \left(1-\frac{1}{l^2}\right)
\prod\limits_{l\dag r} \frac{l(l^2-l-1)}{(l-1)(l^2-1)}.
\end{equation}
Our first result is\\ 

{\bf Theorem 1:}  \begin{it}
Let $r$ be a fixed integer and $A,B\ge 1$. Then, for every 
$c>0$, we have
\begin{eqnarray*}
& &\frac{1}{4AB}\sum\limits_{|a|\le A} \sum\limits_{|b|\le B}\pi_{E(a,b)}^r\\
&=&
C_r\pi_{1/2}(x)+O\left(\left(\frac{1}{A}+\frac{1}{B}\right)x\log x
+\frac{x^{5/4}\log^3 x}{\sqrt{AB}}+
\frac{\sqrt{x}}{\log^c x}\right),
\end{eqnarray*}
where the implied $O$-constant depends only on $c$ and $r$.\end{it}\\

David and Pappalardi \cite{DaP} obtained the above result with 
$(1/A+1/B)x^{3/2}$ in place of 
$(1/A+1/B)x\log x$ and 
$x^{5/2}/(AB)$ in place of $x^{5/4}\log^3 x/\sqrt{AB}$ in the 
$O$-term. 

From Theorem 1, we immediately obtain the following Lang-Trotter type estimate
on average.\\

{\bf Theorem 2:} \begin{it} 
Let $\varepsilon>0$. If $A,B>x^{1/2+\varepsilon}$ and
$AB>x^{3/2+\varepsilon}$, we have as $x\rightarrow\infty$,
\begin{equation} \label{Cor}
\frac{1}{4AB}\sum\limits_{|a|\le A} \sum\limits_{|b|\le B} \pi_{E(a,b)}^r\sim
C_r\frac{\sqrt{x}}{\log x}.
\end{equation}\end{it}

In \cite{DaP}, (\ref{Cor}) was proved under the stronger condition 
$A,B> x^{1+\varepsilon}$. 

David and Pappalardi asked if (\ref{Cor}) is consistent with the Lang-Trotter
conjecture in the sense that
\begin{equation}\label{LT}
\frac{1}{4AB}\sum\limits_{|a|\le A}\sum\limits_{|b|\le B} C_{E(a,b),r} \sim C_r
\end{equation}
as $A,B\rightarrow \infty$. N. Jones \cite{Jon} 
proved that this average estimate 
holds if the summation is restricted to
$a,b$ such that $E(a,b)$ is a Serre curve. An elliptic curve is called a 
Serre curve if $\phi_E($Gal$(\overline{\ratq}/\ratq))$ is an index two subgroup
in GL$_2(\hat{\intz})$, where $\phi_E$ : Gal$(\overline{\ratq}/\ratq)
\rightarrow$ GL$_2(\hat{\intz})$ denotes the Galois representation associated
to $E$. 
By a result of Serre \cite{Ser}, $\phi_E$ is never surjective, so 
in other words, $E$ is a Serre curve if its Galois representation has 
``image as large as possible''.   
Moreover, extending a result of W.D. Duke \cite{Duk}, Jones proved  
that, according to height, 
almost all elliptic curves over $\ratq$ are Serre curves. This gives some
evidence that (\ref{LT}) really holds. 

Furthermore, David and Pappalardi proved that 
$$
\pi_{E(a,b)}^r(x)\sim C_r\sqrt{x}/\log x
$$
holds for ``almost all'' curves $E(a,b)$ with $|a|\le A$ and $|b|\le B$
if $A,B>x^{2+\varepsilon}$ (Theorem 1.3. in \cite{DaP}). 
Here we show that this ``almost-all'' result holds for considerably
smaller $A,B$-ranges.\\ 

{\bf Theorem 3:}\begin{it} Let $\varepsilon>0$ and fix $c>0$. If 
$A,B>x^{1+\varepsilon}$ and 
$x^{3+\varepsilon}<AB<\exp(\exp(\sqrt{x}/\log^c x))$, then
for all $d>2c$ and for all elliptic curves $E(a,b)$ with $|a|\le A$ and
$|b|\le B$ with at most $O(AB/\log^d z)$ exceptions, we have the 
inequality
$$
|\pi_{E(a,b)}^r(x)-C_r\pi_{1/2}(x)|\ll \frac{\sqrt{x}}{\log^c x}.
$$\end{it}

We shall establish the following more general estimate from which Theorem 3
can be derived by the Tur\'an normal order method ({\it c.f.} \cite{DaP}).\\

{\bf Theorem 4:} \begin{it} 
Let $\varepsilon>0$. If $A,B>x^{1/2+\varepsilon}$ and
$AB>x^{3/2+\varepsilon}$, then for every $c>0$, we have  
\begin{eqnarray} \label{T3} 
& &\frac{1}{4AB}\sum\limits_{|a|\le A}\sum\limits_{|b|\le B}
\left| \pi_{E(a,b)}^r(x)-C_r\pi_{1/2}(x)\right|^2\\
&=&
O\left(\left(\frac{1}{A}+\frac{1}{B}\right)x^{2}
+\frac{x^{5/2}\log^3 x }{\sqrt{AB}}+
\frac{x}{\log^c x}+x^{1/2}\log\log (10AB)\right),\nonumber
\end{eqnarray}
where the implied $O$-constant depends only on $c$ and $r$.\end{it}\\

\section{The work of David-Pappalardi}
The following observations are the starting point of David-Pappalardi's work
in \cite{DaP}.\\ 

{\bf Lemma 1:} \begin{it} 
For $r\le 2\sqrt{p}$, the number of $\fie_p$-isomorphism classes of 
elliptic curves over $\fie_p$ with $p+1-r$ points is the total number of 
ideal classes of
the ring $\intz[(D+\sqrt{D})/2]$, where $D=r^2-4p$ is a negative integer which
is congruent to $0$ or $1$ modulo $4$. This number is the Kronecker class 
number $H(r^2-4p)$.\end{it}\\ 

In the following, we set $H_{r,p}=H(r^2-4p)$.\\

{\bf Lemma 2:} \begin{it} 
Suppose that $p\not=2,3$. Then any elliptic curve over $\fie_p$ has a model
$$
E\ :\ Y^2=X^3+aX+b
$$
with $a,b\in \fie_p$. The elliptic curves $E'(a',b')$ over $p$, which are
$\fie_p$-isomorphic to $E$, are given by all the choices
$$
a'=\mu^4a \ \ \ \ \mbox{\rm and } \ \ \ b'=\mu^6b
$$
with $\mu\in \fie_p^*$. The number of such $E'$ is
\begin{eqnarray*}
(p-1)/6, & &  \mbox{\rm if } a=0\ \mbox{\rm and } p\equiv 1\ \mbox{\rm mod } 
3;\\
(p-1)/4, & &  \mbox{\rm if } b=0\ \mbox{\rm and } p\equiv 1\ \mbox{\rm mod } 
4;\\
(p-1)/2, & &  \mbox{\rm otherwise.}
\end{eqnarray*}\end{it}

The above Lemmas 1 and 2 imply that
the number of curves $E(a,b)$ with $a,b\in \intz$, $0\le a,b<p$ and
$a_p(E(a,b))=r$ is
\begin{equation} \label{NC}
\frac{p H_{r,p}}{2}+O(p).
\end{equation}
Now David and Pappalardi \cite{DaP} write
\begin{eqnarray} \label{goal}
& & 
\frac{1}{4AB}\sum\limits_{|a|\le A} \sum\limits_{|b|\le B} \pi^r_{E(a,b)}(x)\\
&=& \frac{1}{4AB}\sum\limits_{B(r)<p\le x} \sharp\{|a|\le A,\ |b|\le B\ :\ 
a_p(E(a,b))=r\},\nonumber
\end{eqnarray}
where $B(r)=\max\{3,r,r^2/4\}$. 
Using (\ref{NC}), the term on the right-hand side is
\begin{equation} \label{DaPa}
\frac{1}{4AB}\sum\limits_{B(r)<p\le x} \left(\frac{2A}{p}+O(1)\right)
\left(\frac{2B}{p}+O(1)\right)\left(\frac{p H_{r,p}}{2}+O(p)\right).
\end{equation}
This asymptotic estimate was 
used by David and Pappalardi to prove 
their main theorem on the average Frobenius distribution of elliptic 
curves (Theorem 1 in
\cite{DaP}). For the main
term in (\ref{DaPa}) David and Pappalardi proved the following.\\ 

{\bf Lemma 3:} \begin{it} Let $r$ be a fixed integer. Then, for every 
$c>0$, we have
$$
\sum\limits_{B(r)<p\le x} \frac{H_{r,p}}{2p}=C_r\pi_{1/2}(x)+
O\left(\frac{\sqrt{x}}{\log^c x}\right),
$$
where the constant $C_r$ is defined as in (\ref{constant}) and 
the implied $O$-constant depends only on $r$ and $c$.\end{it}\\

In this paper we shall sharpen the error term in (\ref{DaPa}).    

\section{Preliminaries}
We first characterize the elliptic curves lying in a fixed 
$\fie_p$-isomorphism class, where $p$ is a prime $\not=2,3$. 
In the following, for $z\in \intz$ let $\overline{z}$ be the 
reduction of $z$ mod $p$. Furthermore, let $z^{-1}$ be a multiplicative inverse
mod $p$, that is, $z z^{-1}\equiv 1$ mod $p$. \\

{\bf Lemma 4:} \begin{it} Let $a,b,c,d\in \intz$, $p\nmid abcd$ and $E_1$,
$E_2$ be  
elliptic curves over $\fie_p$ given by
$$
E_1\ :\ Y^2=X^3+\overline{a}X+\overline{b}.
$$
and 
$$
E_2\ :\ Y^2=X^3+\overline{c}X+\overline{d}.
$$

{\bf (i)} If $p\equiv 1$ mod $4$, then $E_1$ and $E_2$ are 
$\fie_p$-isomorphic if and only if $ca^{-1}$ is a biquadratic residue mod $p$
and $c^3a^{-3}\equiv d^2b^{-2} \mbox{ mod } p$.\medskip

{\bf (ii):} If $p\equiv 3$ mod $4$, then $E_1$ and $E_2$ are 
$\fie_p$-isomorphic
if and only if $ca^{-1}$ and $db^{-1}$ are quadratic
residues mod $p$ and $c^3a^{-3}\equiv d^2b^{-2} \mbox{ mod } p$.\end{it}\\

{\bf Proof:} By Lemma 2, the curves
$E_1$ and $E_2$ are $\fie_p$-isomorphic if and only if there  
exists an integer $m$ such that $p\nmid m$ and 
\begin{equation} \label{bed}
c\equiv m^4a \mbox{ mod } p \ \ \ \ \mbox{\rm and } \ \ \ d\equiv m^6b
\mbox{ mod } p.
\end{equation}

{\bf (i)} Suppose that $p\equiv 1$ mod $4$. 
If (\ref{bed}) is satisfied, then it follows that 
$ca^{-1}$ is a biquadratic residue mod $p$
and $c^3a^{-3}\equiv m^{12}\equiv d^2b^{-2} \mbox{ mod } p$.

Assume, conversely, that  $ca^{-1}$ is a biquadratic residue mod $p$
and 
\begin{equation} \label{congru}
c^3a^{-3}\equiv d^2b^{-2}\mbox{ mod } p. 
\end{equation}
Since $p\equiv 1$ mod $4$, there exist two solutions $m_1,m_2$ of the 
congruence $c\equiv m^4a$ mod $p$ such that $m_2^2\equiv -m_1^2$ mod $p$, 
and (\ref{congru})
implies that $d^2b^{-2}\equiv m_j^{12}$ mod $p$ for $j=1,2$. 
From this it follows 
that $db^{-1}\equiv m_1^6$ mod $p$ or $db^{-1}\equiv -m_1^6\equiv 
m_2^6$ mod $p$. 
Hence, the system
(\ref{bed}) is soluble for $m$. This completes the proof of {\bf (i)}.$\Box$
\medskip

{\bf (ii)} Suppose that $p\equiv 3$ mod $4$. 
If (\ref{bed}) is satisfied, then it follows that 
$ca^{-1}$ and $db^{-1}$ are quadratic
residues mod $p$ and $c^3a^{-3}\equiv m^{12}\equiv d^2b^{-2} \mbox{ mod } p$.

Assume, conversely, that $ca^{-1}$ and $db^{-1}$ are quadratic
residues mod $p$ and (\ref{congru}) is satisfied. Then,
since $p\equiv 3$ mod $4$, $ca^{-1}$ is also a biqadratic residue. 
Hence, there exists a solution $m$ of the 
congruence $c\equiv m^4a$ mod $p$. Further, (\ref{congru})
implies that $d^2b^{-2}\equiv m^{12}$ mod $p$. 
From this it follows that
that $db^{-1}\equiv m^6$ mod $p$ or $db^{-1}\equiv -m^6$ mod $p$. But $-m^6$ is
a quadratic non-residue mod $p$ since $p\equiv 3$ mod $4$. Thus  
$db^{-1}\not\equiv -m^6$ mod $p$ since $db^{-1}$ is supposed to be a quadratic
residue mod $p$. Hence, we have $db^{-1}\equiv m^6$ mod $p$, and so the system
(\ref{bed}) is soluble for $m$. This completes the proof of {\bf (ii)}.$\Box$
\\

We shall detect elliptic curves lying in a fixed 
$\fie_p$-isomorphism class by
using Dirichlet characters. For the estimation of certain error terms we then 
need 
the following results on character sums.\\

{\bf Lemma 5:}\begin{it} Let $q,N\in \mathbbm{N}$ and 
$(a_n)$ be any sequence of complex numbers. Then 
$$
\sum\limits_{\chi\ \! \mbox{\scriptsize\rm mod }q}  
\left| \sum\limits_{n\le N} a_n\chi(n) \right|^2
=\varphi(q)\sum\limits_{\substack{a=1\\ (a,q)=1}}^{q} \left|
\sum\limits_{\substack{n\le N\\ n\equiv a\ \!
\mbox{\scriptsize \rm  mod }q}} a_n\right|^2,
$$ 
where the outer 
sum on the left-hand side runs over all Dirichlet characters mod $q$.\end{it}
\\

{\bf Proof:} This is a consequence of the orthogonality relations for 
Dirichlet characters. $\Box$\\

{\bf Lemma 6:}\begin{it} Let $q,N\in \mathbbm{N}$, $q\ge 2$. Then 
$$
\sum\limits_{\chi\not=\chi_0} \left| \sum\limits_{n\le N} \chi(n) \right|^4
\ll N^2q\log^6 q,
$$
where the outer sum on the left-hand side runs over all non-principal 
Dirichlet characters
mod $q$. \end{it}\\

{\bf Proof:} This is Lemma 3 in \cite{FIw}. $\Box$\\

{\bf Lemma 7:} \begin{it}
Let $q,N\in \mathbbm{N}$, $q\ge 2$ and $\chi$ be any non-principal character 
mod $q$. Then
$$
\sum\limits_{n\le N} \chi(n)\ll \sqrt{q}\log q.
$$\end{it}

{\bf Proof:} This is the well-known inequality of Polya-Vinogradov. $\Box$\\

Furthermore, we shall need the following estimates for sums over 
$H_{r,p}$.\\

{\bf Lemma 8:} \begin{it} We have
$$
\sum\limits_{B(r)<p\le x} H_{r,p}^{1/2}\ll  x^{5/4}, \ \ \ \ \ \
\sum\limits_{B(r)<p\le x} \frac{H_{r,p}}{\sqrt{p}}\ll x, \ \ \ \ 
\ \ \sum\limits_{B(r)<p\le x} \frac{H_{r,p}}{p}\ll \sqrt{x} 
$$ 
and 
$$
\sum\limits_{B(r)<p\le x} \frac{H_{r,p}}{p^2}\ll 1.
$$
\end{it}

{\bf Proof:} By (26) in \cite{DaP}, we have
\begin{equation}\label{HE}
\sum\limits_{B(r)<p\le x} H_{r,p} \ll x^{3/2}.
\end{equation}
Using the Cauchy-Schwarz inequality, we obtain
$$
\sum\limits_{B(r)<p\le x} H_{r,p}^{1/2}\ll 
x^{1/2}\left(\sum\limits_{B(r)<p\le x} H_{r,p}\right)^{1/2}\ll  x^{5/4}
$$
from (\ref{HE}). The remaining three estimates in Lemma 8 can be derived from 
(\ref{HE}) by partial summation.$\Box$\\

Finally, we shall need the following bound.\\

{\bf Lemma 9:} The number of $\fie_p$-isomorphism classes of elliptic curves 
containing curves  
$$
E\ :\ Y^2=X^3+aX+b
$$
over $\fie_p$  with $a=0$ or $b=0$ is bounded by 10.\\

{\bf Proof:} By  Lemma 2, the number of $\fie_p$-isomorphism classes 
containing curves
$E(0,b)$ with $b\in \fie_p^* $ is $\le 6$, and the number 
of $\fie_p$-isomorphism 
classes 
containing curves $E(a,0)$ with $a\in \fie_p^*$ is $\le 4$. $\Box$ 
   
\section{Proof of Theorem 1}
Let $I_{r,p}$ be the number of $\fie_p$-isomorphism classes of elliptic curves 
$$
E\ :\ Y^2=X^3+cX+d
$$
over $\fie_p$ with $p+1-r$ points
such that $c,d\not=0$. Let $(u_{p,j},v_{p,j})$, $j=1,...,I_{r,p}$ be pairs of
integers such that the curves $E(\overline{u_{p,j}},\overline{v_{p,j}})$ 
form a 
system of representatives of these isomophism classes. 
We now write
\begin{eqnarray*}
& & \sharp\{|a|\le A,\ |b|\le B\ :\ a_p(E(a,b))=r\}\\  &=&
\sharp\{|a|\le A,\ |b|\le B\ :\ p\nmid ab,\  a_p(E(a,b))=r\}+\nonumber\\
& & O\left(\frac{AB}{p}+A+B\right)\nonumber
\end{eqnarray*}
and 
\begin{eqnarray}\label{otto}
& & \sharp\{|a|\le A,\ |b|\le B\ :\ p\dag ab,\  a_p(E(a,b))=r\}\\
&=& \sum\limits_{j=1}^{I_{r,p}} 
\sharp\{|a|\le A,\ |b|\le B\ : E(\overline{a},\overline{b})\cong 
E(\overline{u_{p,j}},\overline{v_{p,j}})\},\nonumber
\end{eqnarray}
where  the symbol $\cong$ stands for ``$\fie_p$-isomorphic''.
We rewrite the term on the right-hand  side of (\ref{otto}) as a character sum.
If $p\equiv 1$ mod $4$, then, by Lemma 4(i) and the character relations, 
this term equals
\begin{equation} \label{character}
\frac{1}{4\varphi(p)}\sum\limits_{j=1}^{I_{r,p}} 
\sum\limits_{|a|\le A} \sum\limits_{|b|\le B} \sum\limits_{k=1}^4 
\left(\frac{au_{p,j}^{-1}}{p}\right)_4^k \
\sum\limits_{\chi\ \! \mbox{\scriptsize mod } p} 
\chi(a^3u_{p,j}^{-3}b^{-2}v_{p,j}^2),
\end{equation}
where $(\cdot/p)_4$ is the biquadratic residue symbol. 
If $p\equiv 3$ mod $4$, then, by Lemma 4(ii) and the character relations, 
the term on the right-hand  side of (\ref{otto}) equals
\begin{eqnarray*} & &
\frac{1}{4\varphi(p)}\sum\limits_{j=1}^{I_{r,p}} 
\sum\limits_{|a|\le A} \sum\limits_{|b|\le B} \left(\chi_0(a)+
\left(\frac{au_{p,j}^{-1}}{p}\right)\right) 
\left(\chi_0(b)+
\left(\frac{bv_{p,j}^{-1}}{p}\right)\right) \\
& & \sum\limits_{\chi\ \! \mbox{\scriptsize mod } p} 
\chi(a^3u_{p,j}^{-3}b^{-2}v_{p,j}^2),
\end{eqnarray*}
where $(\cdot/p)$ is the Legendre symbol and $\chi_0$ is the principal
character.

In the following, we consider only the case $p\equiv 1$ mod $4$. The case
$p\equiv 3$ mod $4$ can be treated in a similar way. 
The expression in (\ref{character}) equals
$$
\frac{1}{4\varphi(p)} \sum\limits_{k=1}^4 
\sum\limits_{\chi\ \! \mbox{\scriptsize mod } p} \sum\limits_{j=1}^{I_{r,p}} 
\left(\frac{u_{p,j}}{p}\right)_4^{-k} \overline{\chi}^3(u_{p,j}) 
\chi^2(v_{p,j}) 
\sum\limits_{|a|\le A}\left(\frac{a}{p}\right)_4^{k}\chi^3(a)
\sum\limits_{|b|\le B}  \overline{\chi}^2(b).
$$
We split this expression into 3 parts $M,E_1,E_2$, where\\ \\
({\bf i}) $M=$ contribution of $k,\chi$ with 
$(\cdot/p)_4^{k} \chi^3=\chi_0$, 
$\chi^2=\chi_0$;\\ \\
({\bf ii}) $E_1=$ contribution of $k,\chi$ with 
$(\cdot/p)_4^{k} \chi^3\not=\chi_0$, 
$\chi^2=\chi_0$ or 

\ $(\cdot/p)_4^{k} \chi^3=\chi_0$, 
$\chi^2\not=\chi_0$;\\ \\
({\bf iii}) $E_2=$ contribution of $k,\chi$ with 
$(\cdot/p)_4^{k} \chi^3\not=\chi_0$, $\chi^2\not=\chi_0$.\\ \\
As one may expect, $M$ shall turn out to be the main term and $E_1$, $E_2$
to be the error terms.\\

{\it Estimation of $M$.}\ \ \ \
The only cases in which 
$(\cdot/p)_4^{k} \chi^3=\chi_0 \mbox{ and } \chi^2=\chi_0$ are 
$k=0$, $\chi=\chi_0$ and $k=2$, $\chi=(\cdot/p)$. Now, by a short calculation,
we obtain
\begin{equation} \label{S1}
M=\frac{ABI_{r,p}}{2p}\left(1+O\left(\frac{1}{p}\right)\right).
\end{equation}
By Lemma 9, we have $H_{r,p}-I_{r,p}\le 10$. Combining this with (\ref{S1}), 
we obtain
$$
M=\frac{ABH_{r,p}}{2p}+O\left(\frac{AB}{p}+\frac{ABH_{r,p}}{p^2}\right).
$$\medskip

{\it Estimation of $E_1$.}\ \ \ \ The number of solutions $(k,\chi)$ with
$k=1,...,4$ of 
$(\cdot/p)_4^{k} \chi^3=\chi_0$ is bounded by 12, and $\chi^2=\chi_0$
has precisely 2 solutions $\chi$.
Thus $E_1$ is the sum of at most $12+4\cdot 2=20$
terms of the form
$$ 
\frac{1}{4\varphi(p)} \sum\limits_{j=1}^{I_{r,p}} 
\overline{\chi_1}(u_{p,j})\overline{\chi_2}(v_{p,j})
\sum\limits_{|a|\le A} \chi_1(a) \sum\limits_{|b|\le B} \chi_2(b),
$$
where exactly one of the characters $\chi_1$, $\chi_2$ is the principal
character $\chi_0$. Therefore, Lemma 7 implies that
$$
E_1\ll \frac{I_{r,p}(A+B)}{\sqrt{p}}\log p.
$$ 
 
{\it Estimation of $E_2$.}\ \ \ \  Given $k\in \intz$ and a character
$\chi_1$ mod $p$, 
the number of solutions $\chi$ of 
$\left(\frac{\cdot}{p}\right)_4^{k} {\chi}^{-3}=\chi_1$ 
is $\le 3$, and the number of solutions $\chi$ of $\chi^2=\chi_1$ is $\le 2$. 
Thus, using the Cauchy-Schwarz inequality, 
we deduce that
\begin{eqnarray} \label{cha}
E_2 &\ll& \frac{1}{p} \sum\limits_{k=1}^4 \left(
\sum\limits_{\chi} \left|\sum\limits_{j=1}^{I_{r,p}} 
\left(\frac{u_{p,j}}{p}\right)_4^{k} \chi(u_{p,j}^{-3}v_{p,j}^2)
\right|^2 \right)^{1/2}\times\\ & & \left(\sum\limits_{\chi\not=\chi_0} 
\left| \sum\limits_{|a|\le A} \chi(a)
\right|^4\right)^{1/4} \left(\sum\limits_{\chi\not=\chi_0} \left|
\sum\limits_{|b|\le B}  \chi(b)\right|^4\right)^{1/4}.\nonumber
\end{eqnarray}  
By Lemma 4(i), the number of $j$'s such that 
$u_{p,j}^{-3}v_{p,j}^2$ lie in a fixed residue class mod $p$ is bounded by 4. 
Using this, Lemma 5 and Lemma 6, the expression on the right-hand side of
(\ref{cha}) is dominated by 
$$
\ll (I_{r,p}AB)^{1/2}\log^3 p.
$$

{\it The final estimate.} Combining all contributions, and using 
$I_{r,p}\le H_{r,p}$, we obtain 
\begin{eqnarray}\label{end}
& & \sharp\{|a|\le A,\ |b|\le B\ :\ a_p(E(a,b))=r\}\\  &=&
\frac{ABH_{r,p}}{2p}+O\left(\frac{AB}{p}+\frac{ABH_{r,p}}{p^2}+A+B+
\left(H_{r,p}AB\right)^{1/2}\log^3 p+\right.\nonumber\\
& & \left. \frac{H_{r,p}(A+B)}{\sqrt{p}}\log p\right)
\nonumber
\end{eqnarray}
The result of Theorem 1 now follows from (\ref{goal}), (\ref{end}), Lemma 3
and Lemma 8. 

\section{Proof of Theorem 4} 
As in \cite{DaP}, we set
$$
\mu:=\frac{1}{4AB}\sum\limits_{|a|\le A}\sum\limits_{|b|\le B} 
\pi_{E(a,b)}^r(x).
$$
Fix any $c>0$. Using Theorem 1 and following the arguments in \cite{DaP}, 
if 
$A,B>x^{1/2+\varepsilon}$ and $AB>x^{3/2+\varepsilon}$, then  
\begin{equation}\label{mu}
\mu = C_r\pi_{1/2}(x)+O\left(\frac{\sqrt{x}}{\log^c x}\right),
\end{equation}
and the left-hand side of (\ref{T3}) is
\begin{eqnarray} \label{start}
& & \\
&\ll& \left|\sum\limits_{|a|\le A}\sum\limits_{|b|\le B} \sharp\{
p,q\le x \ :\ p\not=q,\ a_p(E(a,b))=r=a_q(E(a,b))\}\ -\mu^2\right|+\nonumber\\
& & \mu+\frac{x}{\log^{2c} x}\nonumber,
\end{eqnarray}
where $p,q$ denote primes. Similarly as in the preceeding section, we have
\begin{eqnarray} \label{goal2}
& & \sum\limits_{|a|\le A} \sum\limits_{|b|\le B} \sharp\{
p,q\le x \ :\ p\not=q,\ a_p(E(a,b))=r=a_q(E(a,b))\}\\ &=&
\sum\limits_{\substack{B(r)<p,q\le  x\\ p\not= q}} 
\sharp\{|a|\le A,\ |b|\le B\ :\ 
a_p(E(a,b))=r=a_q(E(a,b))\}.\nonumber
\end{eqnarray}
Using Theorem 1 and $\sharp\{p\ :\ p|ab\}=\omega(|ab|)\ll \log\log(10|ab|)$ 
if $ab\not=0$, 
we deduce
\begin{eqnarray} \label{app}
& & \sum\limits_{\substack{B(r)<p,q\le  x\\ p\not= q}}
\sharp\{|a|\le A,\ |b|\le B\ :\ a_p(E(a,b))=r=a_q(E(a,b))\}\\ \nonumber\\
&=& \sum\limits_{\substack{B(r)<p,q\le  x\\ p\not= q}}
\sharp\{|a|\le A,\ |b|\le B\ :\ p,q\nmid ab,\ a_p(E(a,b))=r=a_q(E(a,b))\}
\nonumber\\ & & +O\left(
\sum\limits_{p\le x} \ \sum\limits_{\substack{|a|\le A,\ \! |b|\le B\\
p|ab}} \pi_{E(a,b)}^r(x)\right)
\nonumber\\ &=&
\sum\limits_{\substack{B(r)<p,q\le  x\\ p\not= q}}
\sharp\{|a|\le A,\ |b|\le B\ :\ p,q\nmid ab,\ a_p(E(a,b))=r=a_q(E(a,b))\}
\nonumber\\ & &
+O\left(ABx^{1/2}\log\log(10AB) +(A+B)x^{3/2}\right).\nonumber
\end{eqnarray}

Now we fix $p,q$ with $p\not=q$. 
In the following, we confine ourselves to the case when $p\equiv q\equiv 1$ 
mod $4$. The remaining cases $pq\equiv -1$ mod $4$ and $p\equiv q\equiv 3$ mod
$4$ can be treated in a similar way. 
Similarly as in the preceeding section, we can express the term
$$
\sharp\{|a|\le A,\ |b|\le B\ :\ p,q\nmid ab,\ a_p(E(a,b))=r=a_q(E(a,b))\}
$$
as a character sum
\begin{eqnarray*} \label{character2}
& & \frac{1}{16\varphi(p)\varphi(q)}
\sum\limits_{i=1}^{I_{r,p}}\sum\limits_{j=1}^{I_{r,q}} 
\sum\limits_{|a|\le A} \sum\limits_{|b|\le B} \sum\limits_{k=1}^4
\left(\frac{au_{p,i}^{-1}}{p}\right)_4^k\
\sum\limits_{\chi\ \! \mbox{\scriptsize mod } p} 
\chi(a^3u_{p,i}^{-3}b^{-2}v_{p,i}^2)\\ & &
\sum\limits_{l=1}^4 
\left(\frac{au_{q,j}^{-1}}{q}\right)_4^l\
\sum\limits_{\chi'\ \! \mbox{\scriptsize mod } q} 
\chi'(a^3u_{q,j}^{-3}b^{-2}v_{q,j}^2).\nonumber
\end{eqnarray*}
This sum equals
\begin{eqnarray} \label{character3}
& & \frac{1}{16\varphi(p)\varphi(q)}\sum\limits_{k=1}^4 \sum\limits_{l=1}^4
\ \sum\limits_{\chi\ \! \mbox{\scriptsize mod } p}\
\sum\limits_{\chi'\ \! \mbox{\scriptsize mod } q} 
\left(\sum\limits_{i=1}^{I_{r,p}}
\left(\frac{u_{p,i}}{p}\right)_4^{-k}
{\overline\chi}^3(u_{p,i})\chi^2(v_{p,i})\right)\times\\ & &
\left(\sum\limits_{j=1}^{I_{r,q}}
\left(\frac{u_{q,j}}{q}\right)_4^{-l}
\overline{\chi'}^3(u_{q,j}){\chi'}^2(v_{q,j})\right)
\left(\sum\limits_{|a|\le A} 
\left(\frac{a}{p}\right)_4^k\left(\frac{a}{q}\right)_4^l 
\left(\chi{\chi'}\right)^3(a)\right)\times\nonumber\\ 
& & \left(\sum\limits_{|b|\le B} 
\left(\overline{\chi\chi'}\right)^2(b)\right).\nonumber
\end{eqnarray}
Let $\chi_0$ be the principal character mod $p$ and 
$\chi'_0$ be the principal character mod $q$. Then $\chi_0\chi'_0$ is
the principal character mod $pq$. As previously, we split the
expression in (\ref{character3}) into 3 parts $M,E_1,E_2$, where\\ \\
({\bf i}) $M=$ contribution of $k,l,\chi,\chi'$ with 

$(\cdot/p)_4^{k} (\cdot/q)_4^{l}(\chi\chi')^3=\chi_0\chi'_0$, 
$(\chi\chi')^2=\chi_0\chi'_0$;\\ \\
({\bf ii}) $E_1=$ contribution of $k,l,\chi,\chi'$ with 

\ $(\cdot/p)_4^{k} (\cdot/q)_4^{l}(\chi\chi')^3\not=\chi_0\chi'_0$, 
$(\chi\chi')^2=\chi_0\chi'_0$ or 

\ $(\cdot/p)_4^{k} (\cdot/q)_4^{l}(\chi\chi')^3=\chi_0\chi'_0$, 
$(\chi\chi')^2\not=\chi_0\chi'_0$;\\ \\
({\bf iii}) $E_2=$ contribution of $k,l,\chi,\chi'$ with 

\ \ $(\cdot/p)_4^{k}(\cdot/q)_4^{l}(\chi\chi')^3\not=\chi_0\chi'_0$, 
$(\chi\chi')^2\not=\chi_0\chi'_0$.
\\ 

{\it Estimation of $M$.}\ \ \ \
The only cases in which 
$(\cdot/p)_4^{k} (\cdot/q)_4^{l}(\chi\chi')^3=\chi_0\chi'_0$, 
$(\chi\chi')^2=\chi_0\chi'_0$ are:\medskip\\
(a) $k=l=0$, $\chi=\chi_0$, $\chi'=\chi'_0$; \medskip\\
(b) $k=l=2$, $\chi=(\cdot/p)$, $\chi'=(\cdot/q)$;\medskip\\
(c) $k=0$, $l=2$, $\chi=\chi_0$, $\chi'=(\cdot/q)$;\medskip\\
(d) $k=2$, $l=0$, $\chi=(\cdot/p)$, $\chi'=\chi_0$.\medskip\\ 
Now, by a short calculation, we obtain
\begin{equation} \label{S2}
M=\frac{ABI_{r,p}I_{r,q}}{4pq}\left(1+O\left(\frac{1}{p}+\frac{1}{q}\right)
\right).
\end{equation}
By Lemma 9, we have $H_{r,p}-I_{r,p}\le 10$ and $H_{r,q}-I_{r,q}\le 10$. 
Combining this with (\ref{S2}), 
we obtain
$$ 
M=\frac{ABH_{r,p}H_{r,q}}{4pq}+O\left(\frac{AB(H_{r,p}+H_{r,q})}{pq}
+ABH_{r,p}H_{r,q}\left(\frac{1}{p^2q}+\frac{1}{pq^2}\right)\right).
$$

{\it Estimation of $E_1$.}\ \ \ \ The number of solutions 
$(k,l,\chi,\chi')$ with
$k,l=1,...,4$ of 
$(\cdot/p)_4^{k} (\cdot/q)_4^{l}(\chi\chi')^3\not=\chi_0\chi'_0$
is bounded by $12^2$, and $(\chi\chi')^2=\chi_0\chi_0'$
has precisely 4 solutions $(\chi,\chi')$.
Thus $E_1$ is the sum of at most $144+16\cdot 4=228$
terms of the form
$$ 
\frac{1}{16\varphi(p)\varphi(q)} \sum\limits_{|a|\le A} \chi_1(a) 
\sum\limits_{|b|\le B} \chi_2(b)\sum\limits_{i=1}^{I_{r,p}} 
\chi_3(u_{p,i})\chi_4(v_{p,i}) \sum\limits_{j=1}^{I_{r,q}}
\chi'_3(u_{q,j})\chi'_4(v_{q,j}),
$$
where $\chi_1$, $\chi_2$ are characters mod $pq$ such that exactly one of them
is the principal character, $\chi_3$, $\chi_4$ are characters mod $p$, and
$\chi'_3$, $\chi'_4$ are characters mod $q$. Here the characters 
$\chi_{3,4}$, $\chi'_{3,4}$
depend on the characters $\chi_{1,2}$. Now Lemma 7 implies that
$$
E_1\ll \frac{I_{r,p}I_{r,q}(A+B)}{\sqrt{pq}}\log pq.
$$ 
 
{\it Estimation of $E_2$.}\ \ \ \  Given $k,l\in \intz$ and a character
$\chi_1$ mod $pq$, 
the number of characters $\chi$ mod $pq$ such that  
$(\cdot/p)_4^{k} (\cdot/q)_4^{l}(\chi\chi')^3=\chi_1$  
is $\le 9$, and the number of $\chi$ mod $pq$ such that 
$\chi^2=\chi_1$ is $\le 4$. 
Thus, using the Cauchy-Schwarz inequality, 
we deduce that
\begin{eqnarray} \label{cha2}
& &\\
E_2 &\ll& \frac{1}{pq} \sum\limits_{k=1}^4 \sum\limits_{l=1}^4 \left(
\sum\limits_{\chi} 
\left|\sum\limits_{i=1}^{I_{r,p}} 
\left(\frac{u_{p,i}}{p}\right)_4^{k} \chi(u_{p,i}^{-3}v_{p,i}^2)
\right|^2 \right)^{1/2}\times \nonumber\\ & &
\left(\sum\limits_{\chi'} \left|\sum\limits_{j=1}^{I_{r,q}} 
\left(\frac{u_{q,j}}{q}\right)_4^{l} \chi'(u_{q,j}^{-3}v_{q,j}^2)
\right|^2 \right)^{1/2}
\left(\sum\limits_{\chi_1\not= \chi_0\chi'_0}
\left| \sum\limits_{|a|\le A} \chi_1(a)
\right|^4\right)^{1/4} \nonumber\\ 
& & \left(\sum\limits_{\chi_2\not=\chi_0\chi'_0} 
\left|\sum\limits_{|b|\le B}  \chi(b)\right|^4\right)^{1/4},\nonumber
\end{eqnarray}  
where $\chi$ runs over all characters mod $p$, $\chi'$ runs over all
characters mod $q$, and $\chi_1,\chi_2$ run over all non-principal
characters mod $pq$. 

By Lemma 4(i), the number of $i$'s such that 
$u_{p,i}^{-3}v_{p,i}^2$ lie in a fixed residue class mod $p$ is bounded by 4.
The same is true for the number of $j's$ such that 
$u_{q,j}^{-3}v_{q,j}^2$ lie in a fixed residue class mod $q$.  
Using this, Lemma 5 and Lemma 6, the expression on the right-hand side of
(\ref{cha2}) is dominated by 
$$
\ll (I_{r,p}I_{r,q}AB)^{1/2}\log^3 pq.
$$

{\it The final estimate.} Combining all contributions, and using 
$I_{r,p}\le H_{r,p}$, we obtain 
\begin{eqnarray}\label{end2}
& & \sharp\{|a|\le A,\ |b|\le B\ :\ p,q\nmid ab,
\ a_p(E(a,b))=r=a_q(E(a,b))\}\\  &=&
\frac{ABH_{r,p}H_{r,q}}{4pq}+O\left(\frac{AB(H_{r,p}+H_{r,q})}{pq}
+ABH_{r,p}H_{r,q}\left(\frac{1}{p^2q}+\frac{1}{pq^2}\right)\right. 
\nonumber\\ & & \left. +
(H_{r,p}H_{r,q}AB)^{1/2}\log^3 pq+
\frac{H_{r,p}H_{r,q}(A+B)}{\sqrt{pq}}\log pq\right).
\nonumber
\end{eqnarray}
We have proved this estimate only for distinct primes $p,q$ with $p\equiv q
\equiv 1$ mod $4$, but the same estimate can be proved for
$pq\equiv -1$ mod $4$ and $p\equiv q\equiv 3$ mod $4$ in a similar way. 
Now, from (\ref{app}), (\ref{end2}), Lemma 3 and Lemma 8, we obtain 
\begin{eqnarray}\label{end3}
& & \\
& & \frac{1}{4AB}\sum\limits_{\substack{B(r)<p,q\le  x\\ p\not= q}} 
\sharp\{|a|\le A,\ |b|\le B\ :\ a_p(E(a,b))=r=a_q(E(a,b))\}\nonumber\\ 
&=& (C_r\pi_{1/2}(x))^2 +O\Bigg(\sum\limits_{B(r)<p\le x}
\frac{H_{r,p}^2}{p^2}+\frac{x}{\log^c x}+x^{1/2}\log\log(10AB)
\nonumber\\ 
& & + \frac{x^{5/2}}{\sqrt{AB}}\log^3 x+
\left(\frac{1}{A}+\frac{1}{B}\right)x^{2}\Bigg).
\nonumber
\end{eqnarray}
From (23) in \cite{DaP} and $h(d)\ll\sqrt{|d|}$, we obtain
$H_{r,p}\ll p^{1/2+\varepsilon}$ which implies that 
\begin{equation} \label{triv}
\sum\limits_{B(r)<p\le x} \frac{H_{r,p}^2}{p^2} \ll x^{\varepsilon}.
\end{equation}
The result of Theorem 4 now follows from (\ref{mu}), (\ref{start}),
(\ref{goal2}), (\ref{end3}) and (\ref{triv}). \\ \\ \\
{\bf Acknowledgement.} This paper was written when the author held a
postdoctoral fellowship at the Queen's University in Kingston, 
Canada. The author wishes to thank this institution for financial support. He
would further like to thank 
Prof. Ram Murty for his useful comments.\\

$ $\\
$ $\\
{\bf Address of the author:}\medskip\\ 
Stephan Baier\\ 
Queen's University\\ 
Jeffery Hall\\ 
University Ave\\
Kingston, ON K7L3N6 Canada \medskip\\ 
e-mail: sbaier@mast.queensu.ca
\end{document}